\documentclass[12pt]{article}

\usepackage[geometry]{ifsym}
\usepackage{verbatim}
\usepackage{fancyhdr}
\usepackage{graphicx}
\usepackage{mathrsfs}
\usepackage{amssymb}
\usepackage{pst-poly}

\newtheorem{thm}{Theorem}[section]

\newenvironment{pf}[1][Proof]{\noindent\textbf{#1.} }{\hfill\rule{1mm}{2mm}}

\makeatletter \@addtoreset{equation}{section} \makeatother

\begin{document}

\title{The Algorithmic Complexity of Bondage and Reinforcement Problems in bipartite graphs}
\author{Fu-Tao Hu$^a$, \quad Moo Young Sohn$^b$\footnote{\ Corresponding  author, E-mail address: mysohn@changwon.ac.kr} \\
{\small $^a$School of Mathematical Sciences, Anhui University, Hefei, 230601, P.R. China}\\
{\small Email: hufu@mail.ustc.edu.cn}\\
{\small $^b$Mathematics, Changwon National University, Changwon, 641-773, Republic of Korea}
}

\date{}
\maketitle

\begin{abstract}

Let $G=(V,E)$ be a graph. A subset $D\subseteq V$ is a dominating
set if every vertex not in $D$ is adjacent to a vertex in $D$.
The domination number of $G$, {\red denoted by $\gamma(G)$}, is the smallest cardinality of a dominating set of $G$.
The bondage number of a nonempty graph $G$ is the smallest number
of edges whose removal from $G$ results in a graph with domination
number larger than $\gamma(G)$. The reinforcement number of $G$ is
the smallest number of edges whose addition to $G$
results in a graph with smaller domination number than $\gamma(G)$.
In 2012, Hu and Xu proved that the decision problems for the bondage,
the total bondage, the reinforcement and the total reinforcement
numbers are all NP-hard in general graphs. In this paper, we improve
these results to bipartite graphs.


\end{abstract}

\section{Introduction}

For terminology and notation on graph theory not given here,
the reader is referred to Xu~\cite{xu03}. Let $G=(V,E)$ be a finite,
undirected and simple graph, where $V=V(G)$ is the {\red vertex set} and
$E=E(G)$ is the {\red edge set} of $G$. For a vertex $x\in V(G)$, let
$N_G(x)=\{y: xy\in E(G)\}$ be the {\it open set of neighbors} of $x$ and
$N_G[x]=N_G(x)\cup \{x\}$ be the {\it closed set of neighbors} of $x$.

A subset $D\subseteq V$ is a {\it dominating set} of $G$ if every
vertex in $V-D$ has at least one neighbor in $D$. The {\it
domination number} of $G$, denoted by $\gamma(G)$, is the minimum
cardinality among all dominating sets of $G$. A dominating set $D$ is
called a {\it $\gamma$-set} of $G$ if $|D|=\gamma(G)$.
The domination is an important and {\red classic notion} {\red that has become}
one of the most widely researched topics in graph theory and also is used
to study property of networks frequently. A thorough study of domination appears
in the books~\cite{hh98a, hh98b} by Haynes, Hedetniemi, and Slater.
{\red Among various problems related to} the domination
number, some focus on graph alterations and their effects on the domination number.
Here, we are concerned with {\red two particular graph modifications}, the removal and addition
of edges from a graph. The {\it bondage number} of $G$, denoted by $b(G)$,
is the minimum number of edges whose removal from $G$ results in a graph
{\red with a domination number larger than the one of $G$}. The {\it reinforcement number} of $G$,
denoted by $r(G)$, is the smallest number of edges whose addition to
$G$ results in a graph {\red with a domination number smaller than the one of $G$}.
The bondage
number and the reinforcement number were introduced by Fink et
at.~\cite{fjkr90} and Kok, Mynhardt~\cite{km90}, respectively, in
1990. The reinforcement number for digraphs has been {\red studied} by
Huang, Wang and Xu~\cite{hwx09}. The bondage number and
the reinforcement number are two important parameters for measuring
the vulnerability and stability of the network domination under link
failure and link addition.
Recently, Xu~\cite{xu13} gave a review article on bondage numbers in 2013.

A dominating set $D$ of a graph $G$ without isolated vertices is called a
{\it total dominating set} if every vertex in $D$ is also adjacent to another vertex in $D$.
The {\it total domination number} {\red of $G$}, denoted by $\gamma_t(G)$, is the minimum
cardinality among all total dominating sets of $G$. In this paper,
we use the symbol $D_t$ to denote a total dominating set. A total dominating
set $D_t$ is called a {\it $\gamma_t$-set} of $G$ if $|D_t|=\gamma_t(G)$.
The total domination was introduced by Cockayne et al.~\cite{cdh80}.
Total domination in graphs has been extensively studied in the literature.
In 2009, Henning~\cite{h09} surveyed the recent results on total domination in graphs.
The {\it total bondage number} of $G$ without isolated vertices, denoted by $b_t(G)$,
is the minimum number of edges {\red whose removal from $G$ results in a graph with a
total domination number larger than the one of $G$}. The {\it total
reinforcement number} of $G$ without isolated vertices, denoted by $r_t(G)$,
is the smallest number of edges {\red whose addition from $G$ results in a graph with a
total domination number smaller than the one of $G$}. The total bondage number of a graph
was first studied by Kulli and Patwari~\cite{kp91} and further studied
by Sridharan, Elias, Subramanian~\cite{ses07a}, Huang and Xu~\cite{hx07a}. The total
reinforcement number of a graph was first studied by Sridharan,
Elias, Subramanian~\cite{ses07b} and further studied by Henning, Rad
and Raczek~\cite{hrr11}.

For a graph parameter, knowing whether or not there exists a {\red polynomial-time}
algorithm to compute its exact value is the essential problem.
{\red If the decision problem corresponding to the computation of this parameter} is NP-hard or NP-complete,
then {\red polynomial-time} algorithms for this parameter do not exist unless $NP=P$.
{\red The problem of determining the} domination number has been proved NP-complete
for chordal bipartite graphs~\cite{mb87}. {\red For the total domination number, the problem has
been proved NP-complete for bipartite graphs}~\cite{plh83}.
There are {\red many other complexity results} for variations of domination,
these results can be found in the two books~\cite{c98,hh98b} and
the survey~\cite{h09}.

As regards the bondage problem, Hattingh et al.~\cite{hp08} showed
that the restrained bondage problem is NP-complete even for
bipartite graphs. Hu and Xu~\cite{hx12} have showed that
the bondage, the total bondage, the reinforcement and the total
reinforcement numbers are all NP-hard for general graphs.
{\red We know that even if a problem is known} to be NP-hard or NP-complete,
it may be possible to find a {\red polynomial-time} algorithm for a restricted
set of instances from a particular application. {\red The bondage number and
reinforcement number in graphs are very interesting research problems in
graph theory. There are many results about the bondage number and
reinforcement number in bipartite graphs.} Many famous networks
are bipartite graphs, such as hypercube graphs, partial cube,
grid graphs, median graphs and so on. {\red If we proved these
decision problems for the bondage and the reinforcement are
all NP-hard, then the studies on the bondage number and
reinforcement number in bipartite graphs are more meaningful and
we can directly deduce the decision problems for the bondage and the reinforcement are
both NP-hard in general graphs.} So we should be concerned about
the algorithmic complexity of the bondage and reinforcement {\red problems} in bipartite graphs.

In this paper, we will show that the decision problems for the
bondage, the total bondage, the reinforcement and the total
reinforcement numbers are all NP-hard  even for
bipartite graphs. In other words, there are not
{\red polynomial-time} algorithms to compute these parameters unless $P=NP$. The
proofs are in Section 3, Section 4 and Section 5, respectively.

{\red We have considered about whether these four problems
are belong to NP or not. Since the problem of determining the
domination number is NP-complete, and it is not clear that
there is a polynomial algorithm to verify $\gamma(G-B)>\gamma(G)$ (or
$\gamma(G+R)<\gamma(G)$) for any subset $B\subset E(G)$ (or $R\subset \bar{E(G)}$),
these four problems are not obviously seen to be in NP.
We conjecture that they are not in $NP$.
But we can not prove that determining the bondage and the reinforcement are not NP-problems.
This will be our work to study further. In this paper, we only present the results
that these four problems are all NP-hard in bipartite graphs.}

\section{$3$-satisfiability problem}

In {\it Computers and Intractability: A
Guide to the Theory of NP-Completeness}~\cite{gj79}, Garey and Johnson
outline three steps to prove a decision problem to be NP-hard. We
follow the three steps for proving our four decision problems to be NP-hard.
We prove our results by describing a
polynomial transformation from the known NP-complete problem:
$3$-satisfiability problem. To state the $3$-satisfiability problem,
in this section, we recall some terms.

Let $U$ be a set of Boolean variables. A {\it truth assignment} for
$U$ is a mapping $t: U\to\{T,F\}$. If $t(u)=T$, then $u$ is said to
be ``\,true" under $t$; if $t(u)=F$, then $u$ is said to
be``\,false" under $t$. If $u$ is a variable in $U$, then $u$ and
$\bar{u}$ are {\it literals} over $U$. The literal $u$ is true under
$t$ if and only if the variable $u$ is true under $t$; the literal
$\bar{u}$ is true if and only if the variable $u$ is false.

A {\it clause} over $U$ is a set of literals over $U$. It represents
the disjunction of these literals and is {\it satisfied} by a truth
assignment if and only if at least one of its members is true under
that assignment. A collection $\mathscr C$ of clauses over $U$ is
{\it satisfiable} if and only if there exists some truth assignment
for $U$ that simultaneously satisfies all the clauses in $\mathscr
C$. Such a truth assignment is called a {\it satisfying truth
assignment} for $\mathscr C$. The $3$-satisfiability problem is
specified as follows.

\begin{center}
\begin{minipage}{130mm}
\setlength{\baselineskip}{24pt}

\vskip6pt\noindent {\bf $3$-satisfiability problem (3SAT)}:

\noindent {\bf Instance:}\ {\it A
collection $\mathscr{C}=\{C_1,C_2,\ldots,C_m\}$ of clauses over a
finite set $U$ of variables such that $|C_j| =3$ for $j=1,
2,\ldots,m$.}

\noindent {\bf Question:}\ {\it Is there a truth assignment for $U$
that satisfies all the clauses in $\mathscr{C}$?}

\end{minipage}
\end{center}

\begin{thm} \textnormal{(Theorem 3.1 in~\cite{gj79})}
The $3$-satisfiability problem is NP-complete.
\end{thm}

\section{NP-hardness of bondage}

In this section, we will show that the problem determining the
bondage {\red number in bipartite} graphs is NP-hard. We first state the
problem as the following decision problem.

\begin{center}
\begin{minipage}{130mm}
\setlength{\baselineskip}{24pt}

\vskip6pt\noindent {\bf Bondage problem:}

\noindent {\bf Instance:}\ {\it A nonempty graph $G$ and a positive
integer $k$.}

\noindent {\bf Question:}\ {\it Is $b(G)\le k$?}

\end{minipage}
\end{center}

\vskip6pt\begin{thm}
The bondage problem is NP-hard even when
restricted to bipartite graphs and $k=1$.
\end{thm}

\begin{pf}
Let $U=\{u_1,u_2,\ldots,u_n\}$ and $\mathscr{C}=\{C_1,C_2,
\ldots,C_m\}$ be an arbitrary instance of 3SAT.
A graph $G$ will be constructed from the instance of 3SAT,
such that $\mathscr{C}$ is satisfiable if and only if $b(G)=1$.
Such a graph $G$ can be constructed as follows.

For each variable $u_i\in U$, create a cycle $H_i=(u_i,v_i,\bar{u}_i,r_i,q_i,p_i,u_i)$.
Create a single vertex $c_j$ for each $C_j=\{x_j,y_j,z_j\}\in \mathscr{C}$ and
add {\red the set $E_j=\{c_jx_j, c_jy_j, c_jz_j\}$ to the edge set}. Finally,
add a path $P=s_1s_2s_3$, {\red and} join $s_1$ and $s_3$ to each vertex $c_j$
with $1\le j\le m$.

Figure~\ref{f1} illustrates this {\red construction} when
$U=\{u_1,u_2,u_3,u_4\}$ and $\mathscr{C}=\{C_1,C_2,C_3\}$, where
$C_1=\{u_1,u_2,\bar{u}_3\}, C_2=\{\bar{u}_1,u_2,u_4\},
C_3=\{\bar{u}_2,u_3,u_4\}$.

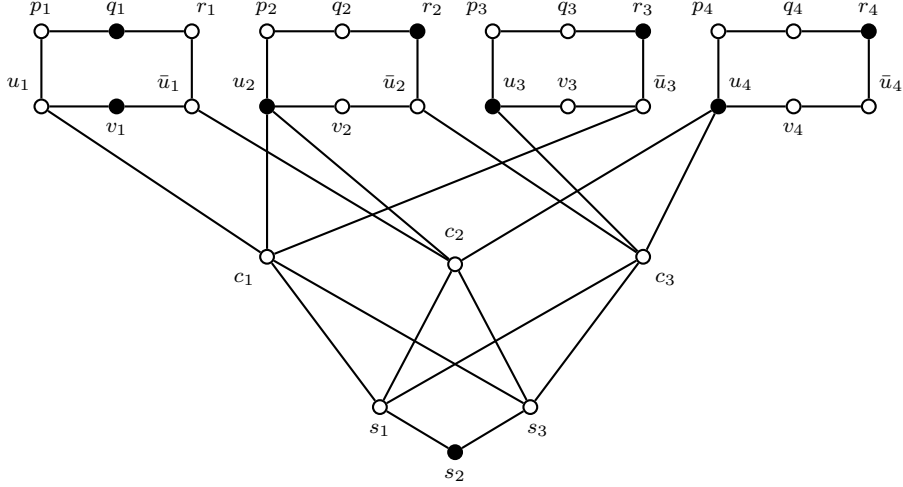
\begin{figure}[ht]
\begin{center}
\begin{pspicture}(-6,-1.1)(6,6.7)

\cnode*(0,-.6){3pt}{s2}\rput(0,-.9){\scriptsize $s_2$}
\cnode(-1,0){3pt}{s1}\rput(-1,-.3){\scriptsize $s_1$}
\cnode(1,0){3pt}{s3}\rput(1.1,-.3){\scriptsize $s_3$} \ncline{s2}{s1}
\ncline{s2}{s3}

\cnode(0,1.9){3pt}{c2}\rput(0,2.3){\scriptsize $c_2$}
\ncline{c2}{s1} \ncline{c2}{s3}
\cnode(-2.5,2){3pt}{c1}\rput(-2.8,1.7){\scriptsize $c_1$}
\ncline{c1}{s1} \ncline{c1}{s3}
\cnode(2.5,2){3pt}{c3}\rput(2.8,1.7){\scriptsize $c_3$}
\ncline{c3}{s1} \ncline{c3}{s3}

\cnode(-5.5,4){3pt}{u1}\rput(-5.8,4.3){\scriptsize $u_1$}
\cnode(-3.5,4){3pt}{u1'}\rput(-3.8,4.3){\scriptsize $\bar{u}_1$}
\cnode*(-4.5,4){3pt}{v1}\rput(-4.5,3.7){\scriptsize $v_1$}
\cnode(-5.5,5){3pt}{p1}\rput(-5.5,5.3){\scriptsize $p_1$}
\cnode(-3.5,5){3pt}{r1}\rput(-3.3,5.3){\scriptsize $r_1$}
\cnode*(-4.5,5){3pt}{q1}\rput(-4.5,5.3){\scriptsize $q_1$}
\ncline{u1}{v1} \ncline{u1}{p1} \ncline{v1}{u1'}
\ncline{u1'}{r1} \ncline{p1}{q1} \ncline{q1}{r1}

\cnode*(-2.5,4){3pt}{u2}\rput(-2.8,4.3){\scriptsize $u_2$}
\cnode(-0.5,4){3pt}{u2'}\rput(-0.8,4.3){\scriptsize $\bar{u}_2$}
\cnode(-1.5,4){3pt}{v2}\rput(-1.5,3.7){\scriptsize $v_2$}
\cnode(-2.5,5){3pt}{p2}\rput(-2.5,5.3){\scriptsize $p_2$}
\cnode*(-.5,5){3pt}{r2}\rput(-.3,5.3){\scriptsize $r_2$}
\cnode(-1.5,5){3pt}{q2}\rput(-1.5,5.3){\scriptsize $q_2$}
\ncline{u2}{v2} \ncline{u2}{p2} \ncline{v2}{u2'}
\ncline{u2'}{r2} \ncline{p2}{q2} \ncline{q2}{r2}

\cnode(2.5,4){3pt}{u3'}\rput(2.8,4.3){\scriptsize $\bar{u}_3$}
\cnode*(0.5,4){3pt}{u3}\rput(0.8,4.3){\scriptsize $u_3$}
\cnode(1.5,4){3pt}{v3}\rput(1.5,4.3){\scriptsize $v_3$}
\cnode*(2.5,5){3pt}{r3}\rput(2.5,5.3){\scriptsize $r_3$}
\cnode(.5,5){3pt}{p3}\rput(.3,5.3){\scriptsize $p_3$}
\cnode(1.5,5){3pt}{q3}\rput(1.5,5.3){\scriptsize $q_3$}
\ncline{u3}{v3} \ncline{u3}{p3} \ncline{v3}{u3'}
\ncline{u3'}{r3} \ncline{p3}{q3} \ncline{q3}{r3}

\cnode(5.5,4){3pt}{u4'}\rput(5.8,4.3){\scriptsize $\bar{u}_4$}
\cnode*(3.5,4){3pt}{u4}\rput(3.8,4.3){\scriptsize $u_4$}
\cnode(4.5,4){3pt}{v4}\rput(4.5,3.7){\scriptsize $v_4$}
\cnode*(5.5,5){3pt}{r4}\rput(5.5,5.3){\scriptsize $r_4$}
\cnode(3.5,5){3pt}{p4}\rput(3.3,5.3){\scriptsize $p_4$}
\cnode(4.5,5){3pt}{q4}\rput(4.5,5.3){\scriptsize $q_4$}
\ncline{u4}{v4} \ncline{u4}{p4} \ncline{v4}{u4'}
\ncline{u4'}{r4} \ncline{p4}{q4} \ncline{q4}{r4}

\ncline{c1}{u1} \ncline{c1}{u2} \ncline{c1}{u3'}
\ncline{c2}{u1'} \ncline{c2}{u2} \ncline{c2}{u4}
\ncline{c3}{u2'} \ncline{c3}{u3} \ncline{c3}{u4}
\end{pspicture}
\caption{\label{f1}\footnotesize An instance of the bondage problem.
Here $\gamma=9$, where the set of bold points is a $\gamma$-set.}
\end{center}
\end{figure}

To prove that this is indeed a transformation, it remains to show that
$b(G)=1$ if and only if there is a truth assignment for $U$ that
satisfies all the clauses in $\mathscr{C}$.
This aim can be {\red fulfilled} by proving the following four claims.

\begin{description}

\item [Claim 3.1]
{\it $\gamma(G)\geq 2n+1$. Moreover, if $\gamma(G)=2n+1$, then for any
$\gamma$-set $D$ in $G$, $D\cap V(P)=\{s_2\}$, $|D\cap V(H_i)|=2$ and
$|D\cap \{u_i, \bar{u}_i\}|\le 1$
for each $i=1,2,\ldots,n$, while $c_j\notin D$ for each
$j=1,2,\ldots,m$.}

\begin{pf}
Let $D$ be a $\gamma$-set of $G$. By the construction of $G$, since
$s_2$ can be dominated only by vertices in $V(P)$, which implies
$|D\cap V(P)|\geq 1$; for each $i=1,2,\ldots,n$, it is
easy to see that $|D\cap N_G[v_i]|\ge 1$ and $|D\cap N_G[q_i]|\ge 1$,
{\red this} implies $|D\cap V(H_i)|\geq 2$. It follows that $\gamma(G)=|D|\geq 2n+1$.

Suppose that $\gamma(G)=2n+1$. Then $|D\cap V(P)|=1$ and $|D\cap
V(H_i)|=2$ for each $i=1,2,\ldots,n$.
Consequently, $c_j\notin D$ for each $j=1,2,\ldots,m$.
Since $q_i$ should be dominated by $D$, $|D\cap \{u_i, \bar{u}_i\}|\le 1$.
Since all vertices in $V(P)$ can be dominated only by $D\cap V(P)$,
{\red this} implies $D\cap V(P)=\{s_2\}$.
\end{pf}

\item [Claim 3.2]
{\it $\gamma(G)=2n+1$ if and only if $\mathscr{C}$ is satisfiable.}

\begin{pf}
Suppose that $\gamma(G)=2n+1$ and let $D$ be a $\gamma$-set of $G$.
By Claim 3.1, for each $i=1,2,\ldots,n$, $|D\cap \{u_i, \bar{u}_i\}|\le 1$.
Define a mapping $t: U\to \{T,F\}$ by
 \begin{equation}\label{e3.1}
 t(u_i)=\left\{
 \begin{array}{l}
 T \ \ {\rm if}\ u_i\in D, \\
 F \ \ {\rm otherwise},
\end{array}
 \right.
 \ i=1,2,\ldots,n.
 \end{equation}

Arbitrarily
choose a clause $C_j\in\mathscr{C}$ with $1\le j\le m$.
There exists some $i$ with $1\le i\le n$ such that $c_j$ is dominated by $u_i\in D$ or
$\bar{u}_i\in D$. Suppose without loss of generality that $c_j$ is dominated by $u_i\in D$.
Since $u_i$ is adjacent to $c_j$ in $G$ and $u_i\in D$, it
follows that $t(u_i)=T$ by (\ref{e3.1}), which implies that the
clause $C_j$ is satisfied by $t$. By the arbitrariness of $j$ with $1\le j\le m$,
it shows that $t$ satisfies all the clauses in $\mathscr{C}$, that is,
$\mathscr{C}$ is satisfiable.

Conversely, suppose that $\mathscr{C}$ is satisfiable, and let $t:
U\to \{T,F\}$ be a satisfying truth assignment for $\mathscr{C}$.
Construct a subset $D'\subseteq V(G)$ as follows. If $t(u_i)=T$,
then put the vertex $u_i$ and $r_i$ in $D'$; if $t(u_i)=F$, then put the
vertex $\bar{u}_i$ and $p_i$ in $D'$. Clearly, $|D'|=2n$. Since $t$ is a
satisfying truth assignment for $\mathscr{C}$, for each
$j=1,2,\ldots,m$, {\red at least one of the three literals} in $C_j$ is true under
the assignment $t$. It follows that $c_j$ {\red can be dominated by $D'$}. Thus
$D'\cup \{s_2\}$ is a dominating set of $G$, and so $\gamma(G)\leq
|D'\cup \{s_2\}|=2n+1$. By Claim 3.1, $\gamma(G)\geq 2n+1$, and so
$\gamma(G)=2n+1$.
\end{pf}

\item [Claim 3.3]
{\it $\gamma(G-e)\leq 2n+2$ for any $e\in E(G)$.}

\begin{pf}
For every edge $e$ in any 6-cycle $H_i$, {\red we have} $\gamma(H_i-e)=2$.
Let $G'$ be the{\red subgraph of $G$ induced} by $\{c_1,c_2,\ldots,c_n,s_1,s_2,s_3\}$ of $G$.
For any edge $e'\in E(G')$, $\{s_1,s_3\}$ is a dominating set of $G'-e'$.
Therefore, $\gamma(G-e)\leq 2n+2$ for any $e\in E(G)$.
\end{pf}

\item [Claim 3.4]
{\it $\gamma(G)=2n+1$ if and only if $b(G)=1$.}

\begin{pf}
Assume $\gamma(G)=2n+1$ and consider the edge $e=s_1s_2$. Suppose
$\gamma(G)=\gamma(G-e)$. Let $D'$ be a $\gamma$-set in $G-e$. It is
clear that $D'$ is also a $\gamma$-set of $G$. By Claim 3.1, we have
$c_j\notin D'$ for each $j=1,2,\ldots,m$ and $D'\cap V(P)=\{s_2\}$.
But then $s_1$ can not be dominated by $D'$, a contradiction. Hence,
$\gamma(G)<\gamma(G-e)$, and so $b(G)=1$.

Now, assume $b(G)=1$. By Claim 3.1, we have that $\gamma(G)\geq
2n+1$. Let $e'$ be an edge such that $\gamma(G)<\gamma(G-e')$. By
Claim 3.3, we have that $\gamma(G-e')\leq 2n+2$. Thus, $2n + 1\leq
\gamma(G)< \gamma(G-e')\leq 2n+2$, which yields $\gamma(G)=2n+1$.
\end{pf}

\end{description}

By Claim 3.2 and Claim 3.4, we prove that $b(G)=1$ if and only if
there is a truth assignment for $U$ that satisfies all the clauses
in $\mathscr{C}$. Since the graph $G$ contains $2n+m+3$ vertices and
$6n+5m+2$ edges, this is clearly a polynomial transformation.
\end{pf}

\section{NP-hardness of total bondage}

In this section, we will show that the problem determining the total
bondage {\red number in} bipartite graphs is NP-hard. We first state it as
the following decision problem.

\begin{center}
\begin{minipage}{130mm}
\setlength{\baselineskip}{24pt}

\vskip6pt\noindent {\bf Total bondage problem:}

\noindent {\bf Instance:}\ {\it A nonempty graph $G$ without isolated vertices and a positive
integer $k$.}

\noindent {\bf Question:}\ {\it Is $b_t(G)\le k$?}

\end{minipage}
\end{center}

\vskip6pt\begin{thm}
The total bondage problem is NP-hard
even when restricted to bipartite graphs and $k=1$.
\end{thm}

\begin{pf}
Let $U=\{u_1,u_2,\ldots,u_n\}$ and $\mathscr{C}=\{C_1,C_2,
\ldots,C_m\}$ be an arbitrary instance of 3SAT. We will construct a graph $G$
such that $\mathscr{C}$ is satisfiable if and only if $b_t(G)=1$.
Such a graph $G$ can be constructed as follows.

For each $u_i\in U$,
{\red create} a graph $H_i$ with {\red vertex set} $V(H_i)=\{u_i,\bar{u}_i,v_i,p_i,q_i\}$
and {\red edge set} $E(H_i)=\{u_iv_i,u_iq_i,\bar{u}_iv_i,v_ip_i,p_iq_i,\bar{u}_iq_i\}$.
For each $C_j=\{x_j,y_j,z_j\}\in \mathscr{C}$, associate a single vertex
$c_j$ and add {\red the set $E_j=\{c_jx_j, c_jy_j,c_jz_j\}$ to the edge set}, $1\le
j\le m$. Finally, add a graph $T$ with {\red vertex set}
$V(T)=\{s_1,s_2,s_3,s_4,s_5,s_6\}$ and {\red edge set}
$E(T)=\{s_1s_2,s_1s_4,s_2s_3,s_2s_5,s_3s_4,s_4s_5,s_5s_6\}$, {\red and} join $s_1$ and $s_3$
to each vertex $c_j$, $1\le j\le m$.

Figure~\ref{f2} shows an example of the graph obtained when
$U=\{u_1,u_2,u_3,u_4\}$ and $\mathscr{C}=\{C_1,C_2,C_3\}$, where
$C_1=\{u_1,u_2,\bar{u}_3\}, C_2=\{\bar{u}_1,u_2,u_4\}$ and $C_3=
\{\bar{u}_2,u_3,u_4\}$.

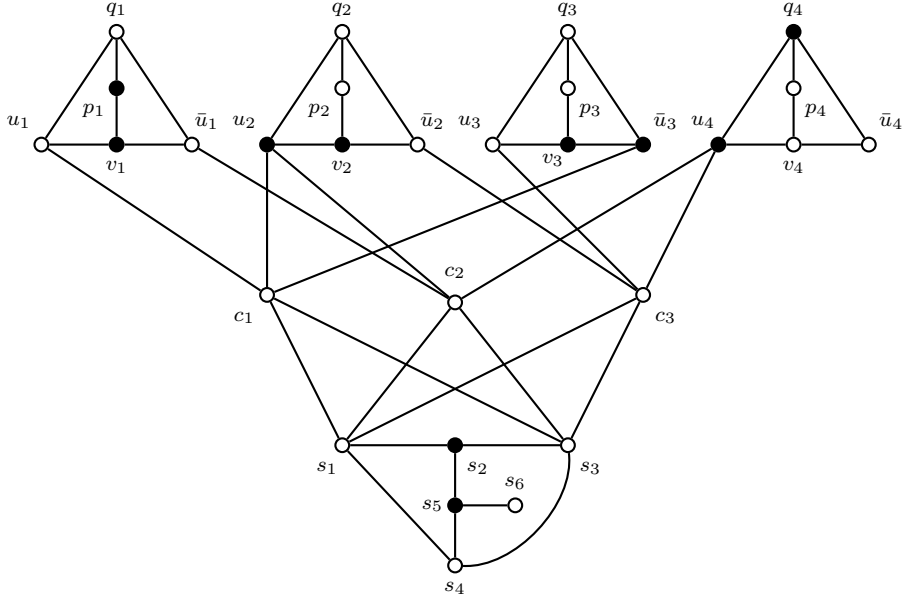
\begin{figure}[ht]
\begin{center}
\begin{pspicture}(-6,-2)(6,7.5)

\cnode*(0,0){3pt}{s2}\rput(.3,-.3){\scriptsize $s_2$}
\cnode(-1.5,0){3pt}{s1}\rput(-1.7,-.3){\scriptsize $s_1$}
\cnode(1.5,0){3pt}{s3}\rput(1.8,-.3){\scriptsize $s_3$} \ncline{s2}{s1} \ncline{s2}{s3}
\cnode(0,-1.6){3pt}{s4}\rput(0,-1.9){\scriptsize $s_4$} \ncline{s4}{s1} \ncarc[arcangle=50]{s3}{s4}
\cnode*(0,-.8){3pt}{s5}\rput(-.3,-.8){\scriptsize $s_5$} \ncline{s5}{s4} \ncline{s5}{s2}
\cnode(.8,-.8){3pt}{s6}\rput(.8,-.5){\scriptsize $s_6$} \ncline{s5}{s6}

\cnode(0,1.9){3pt}{c2}\rput(0,2.3){\scriptsize $c_2$} \ncline{c2}{s1}
\ncline{c2}{s3} \cnode(-2.5,2){3pt}{c1}\rput(-2.8,1.7){\scriptsize $c_1$}
\ncline{c1}{s1} \ncline{c1}{s3}
\cnode(2.5,2){3pt}{c3}\rput(2.8,1.7){\scriptsize $c_3$} \ncline{c3}{s1}
\ncline{c3}{s3}

\cnode(-5.5,4){3pt}{u1}\rput(-5.8,4.3){\scriptsize $u_1$}
\cnode(-3.5,4){3pt}{u1'}\rput(-3.3,4.3){\scriptsize $\bar{u}_1$}
\cnode*(-4.5,4){3pt}{v1}\rput(-4.5,3.7){\scriptsize $v_1$}
\cnode*(-4.5,4.75){3pt}{p1}\rput(-4.8,4.45){\scriptsize $p_1$}
\cnode(-4.5,5.5){3pt}{q1}\rput(-4.5,5.8){\scriptsize $q_1$}
\ncline{u1}{v1} \ncline{u1}{q1} \ncline{v1}{u1'}
\ncline{v1}{p1} \ncline{p1}{q1} \ncline{q1}{u1'}

\cnode*(-2.5,4){3pt}{u2}\rput(-2.8,4.3){\scriptsize $u_2$}
\cnode(-0.5,4){3pt}{u2'}\rput(-0.3,4.3){\scriptsize $\bar{u}_2$}
\cnode*(-1.5,4){3pt}{v2}\rput(-1.5,3.7){\scriptsize $v_2$}
\cnode(-1.5,4.75){3pt}{p2}\rput(-1.8,4.45){\scriptsize $p_2$}
\cnode(-1.5,5.5){3pt}{q2}\rput(-1.5,5.8){\scriptsize $q_2$}
\ncline{u2}{v2} \ncline{u2}{q2} \ncline{v2}{u2'}
\ncline{v2}{p2} \ncline{p2}{q2} \ncline{q2}{u2'}

\cnode*(2.5,4){3pt}{u3'}\rput(2.8,4.3){\scriptsize $\bar{u}_3$}
\cnode(0.5,4){3pt}{u3}\rput(0.2,4.3){\scriptsize $u_3$}
\cnode*(1.5,4){3pt}{v3}\rput(1.3,3.8){\scriptsize $v_3$}
\cnode(1.5,4.75){3pt}{p3}\rput(1.8,4.45){\scriptsize $p_3$}
\cnode(1.5,5.5){3pt}{q3}\rput(1.5,5.8){\scriptsize $q_3$}
\ncline{u3}{v3} \ncline{u3}{q3} \ncline{v3}{u3'}
\ncline{v3}{p3} \ncline{p3}{q3} \ncline{q3}{u3'}

\cnode(5.5,4){3pt}{u4'}\rput(5.8,4.3){\scriptsize $\bar{u}_4$}
\cnode*(3.5,4){3pt}{u4}\rput(3.3,4.3){\scriptsize $u_4$}
\cnode(4.5,4){3pt}{v4}\rput(4.5,3.7){\scriptsize $v_4$}
\cnode(4.5,4.75){3pt}{p4}\rput(4.8,4.45){\scriptsize $p_4$}
\cnode*(4.5,5.5){3pt}{q4}\rput(4.5,5.8){\scriptsize $q_4$}
\ncline{u4}{v4} \ncline{u4}{q4} \ncline{v4}{u4'}
\ncline{v4}{p4} \ncline{p4}{q4} \ncline{q4}{u4'}

\ncline{c1}{u1} \ncline{c1}{u2} \ncline{c1}{u3'}
\ncline{c2}{u1'} \ncline{c2}{u2} \ncline{c2}{u4}
\ncline{c3}{u2'} \ncline{c3}{u3} \ncline{c3}{u4}
\end{pspicture}
\caption{\label{f2}\footnotesize An instance of the total bondage
problem. Here $\gamma_t=10$, where the
set of bold points is a $\gamma_t$-set.}
\end{center}
\end{figure}

It is easy to see that the construction can be accomplished in
polynomial time. All that remains to be shown is that $\mathscr{C}$
is satisfiable if and only if $b_t(G)=1$. This aim can be {\red fulfilled}
by proving the following four claims.

\begin{description}

\item [Claim 4.1]
{\it $\gamma_t(G)\geq 2n+2$. For any $\gamma_t$-set $D_t$ of $G$,
$s_5\in D_t$ and at least one of $v_i$ and $q_i$ {\red belongs to} $D_t$ for each $i=1,2,\ldots,n$. Moreover,
if $\gamma_t(G)=2n+2$, then $D_t\cap V(T)=\{s_2,s_5\}$ or $\{s_4,s_5\}$, $|D_t\cap
V(H_i)|=2$ and $|D_t\cap \{u_i, \bar{u}_i\}|\le 1$ for each $i=1,2,\ldots,n$, while $c_j\notin D_t$ for each
$j=1,2,\dots,m$.}

\begin{pf}
Let $D_t$ be a $\gamma_t$-set of $G$. By the construction of $G$, it
is clear that at least one of $v_i$ and $q_i$ should be in $D_t$ to dominate $p_i$, and
$v_i$ or {\red $q_i$ can} be dominated only by another vertex in $H_i$. It follows
that at least one of $v_i$ and $q_i$ {\red belongs to} $D_t$ and $|D_t\cap V(H_i)|\geq 2$ for each
$i=1,2,\ldots,n$. It is also clear that $s_5$ is certainly in $D_t$
to dominate $s_6$, and $s_5$ can be dominated only by another vertex
in $T$. This fact implies that $s_5\in D_t$ and $|D_t\cap V(T)|\geq
2$. Thus, $\gamma_t(G)=|D_t|\geq 2n+2$.

Suppose that $\gamma_t(G)=2n+2$. Then $|D_t\cap V(H_i)|=2$ for each
$i=1,2,\ldots,n$, and $|D_t\cap V(T)|=2$. Consequently, $c_j\notin
D_t$ for each $j=1,2,\ldots,m$. Since $p_i$ should be dominated by
$D_t$, we have $|D\cap \{u_i, \bar{u}_i\}|\le 1$ for each $i=1,2,\ldots,n$.
{\red Besides}, $s_5$ can be dominated
only by the vertex $s_2$ or $s_4$ in $T$, that is,
at least one of $s_2$ and $s_4$ {\red belongs to} $D_t$. Noting
$|D_t\cap V(T)|=2$, we have $D_t\cap
V(H)=\{s_2,s_5\}~{\red \rm or}~\{s_4,s_5\}$.
\end{pf}

\item [Claim 4.2]
{\it $\gamma_t(G)=2n+2$ if and only if $\mathscr{C}$ is
satisfiable.}

\begin{pf}
Suppose that $\gamma_t(G)=2n+2$ and let $D_t$ be a $\gamma_t$-set of
$G$. By Claim 4.1, $D_t\cap V(T)=\{s_2,s_5\}~{\red \rm or}~\{s_4,s_5\}$ and for each
$i=1,2,\dots,n$, $|D_t\cap \{u_i,\bar{u}_i\}|\le 1$. Define
a mapping $t: U\to \{T,F\}$ by
 \begin{equation}\label{e4.1}
 t(u_i)=\left\{
\begin{array}{ll}
 T \ & {\rm if}\ u_i\in D_t, \\
 F \ & {\rm otherwise},
\end{array}
 \right.
 \ i=1,2,\ldots, n.
 \end{equation}

Arbitrarily choose a clause
$C_j\in\mathscr{C}$. Since the vertex $c_j$ is not
adjacent to any member of $\{s_2, s_4,s_5\}\cup\{v_i,p_i,q_i: 1\le i\le
n\}$, there exists some $i$ with $1\le i\le n$ such that $c_j$ is
dominated by $u_i\in D_t$ or $\bar{u}_i\in D_t$.

Suppose without loss of generality that $c_j$ is dominated by $\bar{u}_i\in D_t$.
Then $\bar{u}_i$ is adjacent to $c_j$ in $G$. Since $\bar{u}_i\in D_t$ and
$|D_t\cap \{u_i,\bar{u}_i\}|\le 1$,
we have $t(\bar{u}_i)=T$ by (\ref{e4.1}), which implies that the clause $C_j$ is
satisfied by $t$. Since the arbitrariness of $j$ with $1\le j\le m$, $\mathscr{C}$ is
satisfiable.

Conversely, suppose that $\mathscr{C}$ is satisfiable, and let $t:
U\to \{T,F\}$ be a satisfying truth assignment for $\mathscr{C}$.
Construct a subset $D'\subseteq V(G)$ as follows. If $t(u_i)=T$,
then put the vertex $u_i$ in $D'$; if $t(u_i)=F$, then put the
vertex $\bar{u}_i$ in $D'$. Clearly, $|D'|=n$. Since $t$ is a
satisfying truth assignment for $\mathscr{C}$, the corresponding vertex $c_j$
in $G$ is adjacent to at least one vertex in $D'$. Let
$D_t'=D'\cup \{s_2,s_5,v_1,\ldots,v_n\}$. Clearly, $D_t'$ is a
total dominating set of $G$ and $|D_t'|=2n+2$. Hence, $\gamma_t(G)\leq |D_t'|=2n+2$. By
Claim 4.1, $\gamma_t(G)\geq 2n+2$. Therefore, $\gamma_t(G)=2n+2$.
\end{pf}

\item [Claim 4.3]
{\it For any $e\in E(G)$, $\gamma_t(G-e)\leq 2n+3$.}

\begin{pf}
It is easy to see that for any edge $e\in E(H_i)$ for each $i=1,2,\ldots,n$,
$\gamma_t(H_i-e)=2$. Let $G'=G-\{H_1,H_2,\ldots,H_n\}$.
For any edge $e'\in E(G')$, {\red it can easily be checked} that $\gamma_t(G')\le 3$.
Thus, for any $e\in E(G)$, $\gamma_t(G-e)\leq 2n+3$.

\end{pf}

\item [Claim 4.4]
{\it $\gamma_t(G)=2n+2$ if and only if $b_t(G)=1$.}

\begin{pf}
Assume $\gamma_t(G)=2n+2$ and take $e=s_2s_5$. Suppose that
$\gamma_t(G-e)=\gamma_t(G)$. Let $D_t'$ be a $\gamma_t$-set of
$G-e$. As $D_t'$ is also a $\gamma_t$-set of $G$, by Claim 4.1 $D_t'\cap
V(H)=\{s_2,s_5\}~{\red \rm or}~\{s_4,s_5\}$, which contradicts {\red the fact that $s_2$
is dominated by $D_t'$} in $G-e$. This contradiction
shows that $\gamma_t(G-e)>\gamma_t(G)$, hence $b_t(G)=1$.

Now, assume $b_t(G)=1$. By Claim 4.1, we have that $\gamma_t(G)\geq
2n+2$. Let $e'$ be an edge such that $\gamma_t(G-e')>\gamma_t(G)$.
By Claim 4.3, we have that $\gamma_t(G-e)\leq 2n+3$. Thus, $2n +
2\leq \gamma_t(G)< \gamma_t(G-e')\leq 2n+3$, which yields
$\gamma_t(G)=2n+2$.
\end{pf}

\end{description}

It follows from Claim 4.2 and Claim 4.4 that $b_t(G)=1$ if and only
if $\mathscr{C}$ is satisfiable. The theorem follows.
\end{pf}

\section{ NP-hardness of reinforcement}

In this section, we will show that {\red the problems of determining the
reinforcement number and total reinforcement number in bipartite graphs are NP-hard}.
We first state them as the following decision problem.

\begin{center}
\begin{minipage}{130mm}
\setlength{\baselineskip}{24pt}

\vskip6pt\noindent {\bf (Total) Reinforcement problem:}

\noindent {\bf Instance:}\ {\it A graph $G$ and a positive integer
$k$.}

\noindent {\bf Question:}\ {\it Is $(r_t(G))\,r(G)\le k$?}

\end{minipage}
\end{center}

\vskip6pt
\begin{thm}\label{thm5.1}
The reinforcement problem is NP-hard
even when restricted to bipartite graphs and $k=1$.
\end{thm}

\begin{pf}
Let $U=\{u_1,u_2,\ldots,u_n\}$ and $\mathscr{C}=\{C_1,C_2,
\ldots,C_m\}$ be an arbitrary instance of 3SAT. We will construct a graph
$G$ such that $\mathscr{C}$ is satisfiable if and only
if $r(G)=1$. Such a graph $G$ can be constructed as follows.

For each $u_i\in U$,
associate a cycle $H_i=(u_i,v_i,\bar{u}_i,r_i,q_i,p_i,u_i)$.
For each $C_j=\{x_j,y_j,z_j\}\in \mathscr{C}$, associate a single vertex $c_j$
and add edges $(c_j,x_j), (c_j,y_j)$ and $(c_j,z_j)$, $1\le j\le m$. Finally, add
a vertex $s$ and join $s$ to every vertex $c_j$.

Figure~\ref{f3} shows an example of the graph obtained when
$U=\{u_1,u_2,u_3,u_4\}$ and $\mathscr{C}=\{C_1,C_2,C_3\}$, where
$C_1=\{u_1,u_2,\bar{u}_3\}, C_2=\{\bar{u}_1,u_2,u_4\},
C_3=\{\bar{u}_2,u_3,u_4\}$.

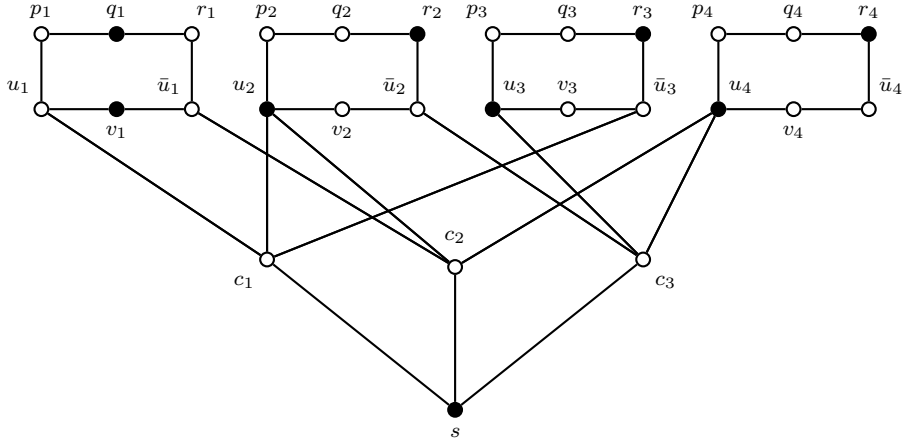
\begin{figure}[ht]
\begin{center}
\begin{pspicture}(-5,-.5)(5,6.7)

\cnode*(0,0){3pt}{s}\rput(0,-.3){\scriptsize $s$}

\cnode(0,1.9){3pt}{c2}\rput(0,2.3){\scriptsize $c_2$} \ncline{c2}{s}
\cnode(-2.5,2){3pt}{c1}\rput(-2.8,1.7){\scriptsize $c_1$} \ncline{c1}{s}
\cnode(2.5,2){3pt}{c3}\rput(2.8,1.7){\scriptsize $c_3$} \ncline{c3}{s}

\cnode(-5.5,4){3pt}{u1}\rput(-5.8,4.3){\scriptsize $u_1$}
\cnode(-3.5,4){3pt}{u1'}\rput(-3.8,4.3){\scriptsize $\bar{u}_1$}
\cnode*(-4.5,4){3pt}{v1}\rput(-4.5,3.7){\scriptsize $v_1$}
\cnode(-5.5,5){3pt}{p1}\rput(-5.5,5.3){\scriptsize $p_1$}
\cnode(-3.5,5){3pt}{r1}\rput(-3.3,5.3){\scriptsize $r_1$}
\cnode*(-4.5,5){3pt}{q1}\rput(-4.5,5.3){\scriptsize $q_1$}
\ncline{u1}{v1} \ncline{u1}{p1} \ncline{v1}{u1'}
\ncline{u1'}{r1} \ncline{p1}{q1} \ncline{q1}{r1}

\cnode*(-2.5,4){3pt}{u2}\rput(-2.8,4.3){\scriptsize $u_2$}
\cnode(-0.5,4){3pt}{u2'}\rput(-0.8,4.3){\scriptsize $\bar{u}_2$}
\cnode(-1.5,4){3pt}{v2}\rput(-1.5,3.7){\scriptsize $v_2$}
\cnode(-2.5,5){3pt}{p2}\rput(-2.5,5.3){\scriptsize $p_2$}
\cnode*(-.5,5){3pt}{r2}\rput(-.3,5.3){\scriptsize $r_2$}
\cnode(-1.5,5){3pt}{q2}\rput(-1.5,5.3){\scriptsize $q_2$}
\ncline{u2}{v2} \ncline{u2}{p2} \ncline{v2}{u2'}
\ncline{u2'}{r2} \ncline{p2}{q2} \ncline{q2}{r2}

\cnode(2.5,4){3pt}{u3'}\rput(2.8,4.3){\scriptsize $\bar{u}_3$}
\cnode*(0.5,4){3pt}{u3}\rput(0.8,4.3){\scriptsize $u_3$}
\cnode(1.5,4){3pt}{v3}\rput(1.5,4.3){\scriptsize $v_3$}
\cnode*(2.5,5){3pt}{r3}\rput(2.5,5.3){\scriptsize $r_3$}
\cnode(.5,5){3pt}{p3}\rput(.3,5.3){\scriptsize $p_3$}
\cnode(1.5,5){3pt}{q3}\rput(1.5,5.3){\scriptsize $q_3$}
\ncline{u3}{v3} \ncline{u3}{p3} \ncline{v3}{u3'}
\ncline{u3'}{r3} \ncline{p3}{q3} \ncline{q3}{r3}

\cnode(5.5,4){3pt}{u4'}\rput(5.8,4.3){\scriptsize $\bar{u}_4$}
\cnode*(3.5,4){3pt}{u4}\rput(3.8,4.3){\scriptsize $u_4$}
\cnode(4.5,4){3pt}{v4}\rput(4.5,3.7){\scriptsize $v_4$}
\cnode*(5.5,5){3pt}{r4}\rput(5.5,5.3){\scriptsize $r_4$}
\cnode(3.5,5){3pt}{p4}\rput(3.3,5.3){\scriptsize $p_4$}
\cnode(4.5,5){3pt}{q4}\rput(4.5,5.3){\scriptsize $q_4$}
\ncline{u4}{v4} \ncline{u4}{p4} \ncline{v4}{u4'}
\ncline{u4'}{r4} \ncline{p4}{q4} \ncline{q4}{r4}

\ncline{c1}{u1} \ncline{c1}{u2} \ncline{c1}{u3'}
\ncline{c2}{u1'} \ncline{c2}{u2} \ncline{c2}{u4}
\ncline{c3}{u2'} \ncline{c3}{u3} \ncline{c3}{u4}

\ncline{c1}{u1} \ncline{c1}{u2} \ncline{c1}{u3'}
\ncline{c2}{u1'} \ncline{c2}{u2} \ncline{c2}{u4}
\ncline{c3}{u2'} \ncline{c3}{u3} \ncline{c3}{u4}
\end{pspicture}
\caption{\label{f3}\footnotesize An instance of the reinforcement
problem. Here $\gamma=5$, where the set
of bold points is a $\gamma$-set.}
\end{center}
\end{figure}

It is easy to see that the construction can be accomplished in
polynomial time. All that remains to be shown is that $\mathscr{C}$
is satisfiable if and only if $r(G)=1$. To this aim, we first prove
the following two claims.

\begin{description}

\item [Claim 5.1.1]
{\it $\gamma(G)=2n+1$.}

\begin{pf}
On the one hand, let $D$ be a $\gamma$-set of $G$, then $\gamma(G)=|D|\geq 2n+1$
since $|D\cap V(H_i)|\geq 2$ and $|D\cap N[s]|\geq 1$. On the other
hand, $D'=\{s,u_1,r_1,u_2,r_2,\ldots,u_n,r_n\}$ is a dominating set of $G$,
which implies that $\gamma(G)\le |D'|=2n+1$. It follows that
$\gamma(G)=2n+1$.
\end{pf}

\item [Claim 5.1.2]
{\it If there exists an edge $e\in E(\bar{G})$ such that
$\gamma(G+e)=2n$, {\red and if $D_e$ denotes a $\gamma$-set of $G + e$}, then
$|D_e\cap V(H_i)|=2$ and $|D_e\cap \{u_i, \bar{u}_i\}|\le 1$ for each $i=1,2,\ldots,n$,
while $s\notin D_e$ and $c_j\notin D_e$ for each $j=1,2,\ldots,m$.}

\begin{pf}
Suppose to the contrary that $|D_e\cap V(H_{i_0})|<2$ for some $i_0$
with $1\le i_0\le n$. Since $\{v_{i_0},p_{i_0},q_{i_0},r_{i_0}\}$ should be
dominated by $D_e$, $D_e\cap V(H_{i_0})=\{q_{i_0}\}$, and then one end-vertex of the edge $e$ should be
$v_{i_0}$ since $D_e$ dominates it via the edge $e$ in $G+e$, and
for every $i\neq i_0$, $|D_e\cap V(T_i)|\geq 2$ since $D_e$
dominates $\{v_i,p_i,q_i,r_i\}$. By the hypotheses, two literals $u_{i_0}$ and
$\bar{u}_{i_0}$ do not simultaneously appear in the same clause in
$\mathscr{C}$, {\red there is no $j$ such that vertex $c_j$ is adjacent to both of them}.
Since $u_{i_0}$ and $\bar{u}_{i_0}$ should be
dominated by $D_e$, there exist two distinct vertices $c_j, c_l\in
D_e$ such that $c_j$ dominates $u_{i_0}$ and $c_l$ dominates
$\bar{u}_{i_0}$. Thus, $|D_e|\geq 2n+1$, a contradiction. Hence,
$|D_e\cap V(H_i)|=2$ for each $i=1,2,\ldots,n$, and $c_j\notin D_e$
for every $j$ since $|D_e|=2n$. Therefore, $s$ should be
dominated by $D_e$ via the edge $e$ in $G+e$. Since $q_i$ should be
dominated by $D_e$, $|D_e\cap \{u_i, \bar{u}_i\}|\le 1$ for each $i=1,2,\ldots,n$.
\end{pf}

\end{description}

We now show that $\mathscr{C}$ is satisfiable if and only if
$r(G)=1$.

Suppose that $\mathscr{C}$ is satisfiable, and let $t: U\to \{T,F\}$
be a satisfying truth assignment for $\mathscr{C}$. We construct a
subset $D'\subseteq V(G)$ as follows. If $t(u_i)=T$ then put the
vertex $u_i$ and $r_i$ in $D'$; if $t(u_i)=F$ then put the vertex
$\bar{u}_i$ and $p_i$ in $D'$. Then $|D'|=2n$. Since $t$ is a satisfying truth assignment
for $\mathscr{C}$, for each $j=1,2,\ldots,m$, at least one of {\red the three}
literals in $C_j$ is true under the assignment $t$. It follows that
the corresponding vertex $c_j$ in $G$ is adjacent to at least one
vertex in $D'$ since $c_j$ is adjacent to each literal in $C_j$ by
the construction of $G$. Without loss of generality let $t(u_1)=T$,
then $D'$ is a dominating set of $G+su_1$, and hence
$\gamma(G+su_1)\leq |D'|=2n$. By Claim 5.1.1, we have
$\gamma(G)=2n+1$. It follows that $\gamma(G+su_1)\leq
2n<2n+1=\gamma(G)$, which implies $r(G)=1$.

Conversely, assume $r(G)=1$. Then there exists an edge $e$ in
$\bar{G}$ such that $\gamma(G+e)=2n$. Let $D_e$ be a $\gamma$-set of
$G+e$. By Claim 5.1.2, {\red $|D_e \{u_i,\bar{u}_i\}|\le 1$} for each
$i=1,2,\ldots,n$, $s\notin D_e$ and $c_j\notin D_e$ for each $j=1,2,\ldots,m$.
Define $t: U\to \{T,F\}$ by
 \begin{equation}\label{e5.1}
 t(u_i)=\left\{
\begin{array}{ll}
 T \ & {\rm if}\ u_i\in D_e, \\
 F \ & {\rm otherwise},
\end{array}
 \right.
 \ i=1,2,\ldots,n.
 \end{equation}

We will show that $t$ is a satisfying truth assignment for
$\mathscr{C}$. It is sufficient to show that every clause in
$\mathscr{C}$ is satisfied by $t$.  

Consider arbitrary clause $C_j\in\mathscr{C}$ with $1\le j\le m$. By
Claim 5.1.2, the corresponding vertex $c_j$ in $G$ is dominated by
$u_i$ or $\bar{u}_i$ in $D_e$ for some $i$. Suppose without loss of generality that $c_j$ is
dominated by $u_i\in D_e$. Then $u_i$ is adjacent to $c_j$ in $G$,
that is, the literal $u_i$ is in the clause $C_j$ by the
construction of $G$. Since $u_i\in D_e$, we have $t(u_i)=T$ by
(\ref{e5.1}), which implies that $C_j$ is satisfied by $t$. The
arbitrariness of $j$ with $1\le j\le m$ shows that {\red every clause}
in $\mathscr{C}$ is satisfied by $t$, that is, $\mathscr{C}$ is
satisfiable.
\end{pf}

\vskip6pt

By using an analogous argument as in the proof of
Theorem~\ref{thm5.1}, we can prove that total reinforcement problem
is also NP-hard even when restricted to bipartite graphs and $k=1$.
Here we give an outline of the proof, the details are omitted.

\vskip6pt
\begin{thm}
The total reinforcement problem is NP-hard
even when restricted to bipartite graphs and $k=1$.
\end{thm}

\begin{pf}
Let $U=\{u_1,u_2,\ldots,u_n\}$ and $\mathscr{C}=\{C_1,C_2,
\ldots,C_m\}$ be an arbitrary instance of 3SAT. We will construct a graph $G$
such that $\mathscr{C}$ is satisfiable if and only if $r_t(G)=1$.
Such a graph $G$ can be constructed as follows.

For each $u_i\in U$, associate a graph $H_i$ with {\red vertex set}
$V(H_i)=\{u_i,\bar{u}_i,v_i,p_i,q_i\}$ and {\red edge set}
$E(H_i)=\{u_iv_i,u_iq_i,\bar{u}_iv_i,v_ip_i,p_iq_i,\bar{u}_iq_i\}$.
For each $C_j=\{x_j,y_j,z_j\}\in \mathscr{C}$, associate a single vertex
$c_j$ and add an {\red edge set} $E_j=\{c_jx_j, c_jy_j,c_jz_j\}$, $1\le
j\le m$. Finally, add a path $P=s_1s_2s_3$ and join $s_1$ to each
vertex $c_j$, $1\le j\le m$.

Figure~\ref{f4} shows an example of the graph obtained when
$U=\{u_1,u_2,u_3,u_4\}$ and $\mathscr{C}=\{C_1,C_2,C_3\}$, where
$C_1=\{u_1,u_2,\bar{u}_3\}, C_2=\{u_1,\bar{u}_2,u_4\}$ and $C_3=
\{\bar{u}_2,\bar{u}_3,u_4\}$.

\begin{figure}[ht]
\begin{center}
\begin{pspicture}(-6,-1)(6,7.5)

\cnode*(0,0){3pt}{s1}\rput(.3,-.2){\scriptsize $s_1$}
\cnode*(0,-.8){3pt}{s2}\rput(-.3,-1){\scriptsize $s_2$} \ncline{s1}{s2}
\cnode(.8,-.8){3pt}{s3}\rput(.8,-.5){\scriptsize $s_3$} \ncline{s2}{s3}

\cnode(0,1.9){3pt}{c2}\rput(.2,2.3){\scriptsize $c_2$} \ncline{c2}{s1}
\cnode(-2.5,2){3pt}{c1}\rput(-2.8,1.7){\scriptsize $c_1$} \ncline{c1}{s1}
\cnode(2.5,2){3pt}{c3}\rput(2.8,1.7){\scriptsize $c_3$} \ncline{c3}{s1}

\cnode*(-5.5,4){3pt}{u1}\rput(-5.8,4.3){\scriptsize $u_1$}  \ncarc[linecolor=red,arcangle=30]{s2}{u1}
\cnode(-3.5,4){3pt}{u1'}\rput(-3.3,4.3){\scriptsize $\bar{u}_1$}
\cnode*(-4.5,4){3pt}{v1}\rput(-4.4,3.8){\scriptsize $v_1$}
\cnode(-4.5,4.75){3pt}{p1}\rput(-4.8,4.45){\scriptsize $p_1$}
\cnode(-4.5,5.5){3pt}{q1}\rput(-4.5,5.8){\scriptsize $q_1$}
\ncline{u1}{v1} \ncline{u1}{q1} \ncline{v1}{u1'}
\ncline{v1}{p1} \ncline{p1}{q1} \ncline{q1}{u1'}

\cnode(-2.5,4){3pt}{u2}\rput(-2.8,4.3){\scriptsize $u_2$}
\cnode(-0.5,4){3pt}{u2'}\rput(-0.3,4.3){\scriptsize $\bar{u}_2$}
\cnode*(-1.5,4){3pt}{v2}\rput(-1.5,3.7){\scriptsize $v_2$}
\cnode*(-1.5,4.75){3pt}{p2}\rput(-1.8,4.45){\scriptsize $p_2$}
\cnode(-1.5,5.5){3pt}{q2}\rput(-1.5,5.8){\scriptsize $q_2$}
\ncline{u2}{v2} \ncline{u2}{q2} \ncline{v2}{u2'}
\ncline{v2}{p2} \ncline{p2}{q2} \ncline{q2}{u2'}

\cnode*(2.5,4){3pt}{u3'}\rput(2.8,4.3){\scriptsize $\bar{u}_3$}
\cnode(0.5,4){3pt}{u3}\rput(0.2,4.3){\scriptsize $u_3$}
\cnode*(1.5,4){3pt}{v3}\rput(1.3,3.8){\scriptsize $v_3$}
\cnode(1.5,4.75){3pt}{p3}\rput(1.8,4.45){\scriptsize $p_3$}
\cnode(1.5,5.5){3pt}{q3}\rput(1.5,5.8){\scriptsize $q_3$}
\ncline{u3}{v3} \ncline{u3}{q3} \ncline{v3}{u3'}
\ncline{v3}{p3} \ncline{p3}{q3} \ncline{q3}{u3'}

\cnode(5.5,4){3pt}{u4'}\rput(5.8,4.3){\scriptsize $\bar{u}_4$}
\cnode*(3.5,4){3pt}{u4}\rput(3.3,4.3){\scriptsize $u_4$}
\cnode(4.5,4){3pt}{v4}\rput(4.5,3.7){\scriptsize $v_4$}
\cnode(4.5,4.75){3pt}{p4}\rput(4.8,4.45){\scriptsize $p_4$}
\cnode*(4.5,5.5){3pt}{q4}\rput(4.5,5.8){\scriptsize $q_4$}
\ncline{u4}{v4} \ncline{u4}{q4} \ncline{v4}{u4'}
\ncline{v4}{p4} \ncline{p4}{q4} \ncline{q4}{u4'}

\ncline{c1}{u1} \ncline{c1}{u2} \ncline{c1}{u3'}
\ncline{c2}{u1} \ncline{c2}{u2'} \ncline{c2}{u4}
\ncline{c3}{u2'} \ncline{c3}{u3'} \ncline{c3}{u4}
\end{pspicture}
\caption{\label{f4}\footnotesize An instance of the total reinforcement
problem. Here $\gamma_t=10$, where the
set of bold points is a $\gamma_t$-set. {\red Add the edge $u_1s_2$ and remove the vertex $s_1$ to decrease the total domination number}.}
\end{center}
\end{figure}

It is easy to see that the construction can be accomplished in
polynomial time. All that remains to be shown is that $\mathscr{C}$
is satisfiable if and only if $r_t(G)=1$.

\begin{description}

\item [Claim 5.2.1]
{\it $\gamma_t(G)=2n+2$.}

\item [Claim 5.2.2]
If there exists an edge $e\in E(\bar{G})$ such that
$\gamma_t(G+e)<2n+2$, and {\red if} $D_e$ be a $\gamma_t$-set of $G+e$,
then $|D_e\cap V(H_i)|=2$ and $|D_e\cap \{u_i, \bar{u}_i\}|\le 1$ for each $i=1,2,\ldots,n$, while
$s_1\notin D_e$ and $c_j\notin D_e$ for each $j=1,2,\ldots,m$.
\end{description}

{\red We now show that $\mathscr{C}$ is satisfiable if and only if
$r_t(G)=1$.

Suppose that $\mathscr{C}$ is satisfiable, and let $t: U\to \{T,F\}$
be a satisfying truth assignment for $\mathscr{C}$. We construct a
subset $D'\subseteq V(G)$ as follows. If $t(u_i)=T$ then put the
vertex $u_i$ in $D'$; if $t(u_i)=F$ then put the vertex
$\bar{u}_i$ in $D'$. Then $|D'|=n$. Let $D_t'=D'\cup \{v_1,v_2,\ldots,v_n,s_2\}$.
Without loss of generality let $u_1\in D_t'$.
We can easily check that $D_t'$ is a total dominating set of $G+s_2u_1$, and hence
$\gamma_t(G+s_2u_1)\leq |D_t'|=2n+1$. By Claim 5.2.1, we have
$\gamma_t(G)=2n+2$. It follows that $r_t(G)=1$.

Conversely, assume $r_t(G)=1$. Then there exists an edge $e$ in
$\bar{G}$ such that $\gamma(G+e)=2n$. Let $D_e$ be a $\gamma_t$-set of
$G+e$. By Claim 5.1.2, {\red $|D_e \{u_i,\bar{u}_i\}|\le 1$} for each
$i=1,2,\ldots,n$, $s_1\notin D_e$ and $c_j\notin D_e$ for each $j=1,2,\ldots,m$.
Define $t: U\to \{T,F\}$ by
\begin{equation}
 t(u_i)=\left\{
\begin{array}{ll}
 T \ & {\rm if}\ u_i\in D_e, \\
 F \ & {\rm otherwise},
\end{array}
 \right.
 \ i=1,2,\ldots,n.
 \end{equation}

Using the same methods as in Theorem~\ref{thm5.1},
we can show that $t$ is a satisfying truth assignment for
$\mathscr{C}$.}
\end{pf}

\section*{Acknowledgments.}
{\red The authors would like to thank the anonymous referees for their
kind comments and helpful suggestions on the original manuscript,
which resulted in this revised version.}

This research was supported by Basic Science Research Program through the National Research Foundation of {\red Korea (NRF)} funded by the Ministry of Education, Science and {\red Technolog (2012R1A1A2005115)}. The first author was supported by the doctoral scientific research startup fund of {\red Anhui  University}.


\begin{thebibliography}{99}

\bibitem{c98}
G.J. Chang, Algorithmic aspects of domination in graphs,
in: D.-Z. Du and P. M. Pardalos, (Ed.), {\it Handbook of Combinatorial
Optimization}, Vol. 3, (1998), 339-405.

\bibitem{cdh80}
E.J. Cockayne, R.M. Dawes, S.T. Hedetniemi, Total domination in
graphs. {\it Networks}, {\bf 10} (1980), 211-219.


\bibitem{fjkr90}
J.F. Fink, M.S. Jacobson, L.F. Kinch, J. Roberts, The bondage
number of a graph. {\it Discrete Math.},  {\bf 86} (1990),
47-57.

\bibitem{gj79}
M.R. Garey, D.S. Johnson, {\it Computers and Intractability: A
Guide to the Theory of NP-Completeness}, Freeman, San Francisco,
1979.

\bibitem{hp08}
J.H. Hattingh, A. R. Plummer, Restrained bondage in graphs. {\it
Discrete Math.}, {\bf 308} (2008), 5446-5453.

\bibitem{h09}
M.A. Henning, A survey of selected recent results on total
domination in graphs. {\it Discrete Math.}, {\bf 309}(1)
(2009), 32-63.

\bibitem{hh98a}
T.W. Haynes, S.T. Hedetniemi, P.J. Slater, {\it Fundamentals of
Domination in Graphs}, Marcel Dekker, New York, 1998.

\bibitem{hh98b}
T.W. Haynes, S.T. Hedetniemi, P.J. Slater, {\it Domination in
Graphs: Advanced Topics}, Marcel Dekker, New York, 1998.

\bibitem{hrr11}
M.A. Henning, N.J. Rad, J. Raczek, A note on
total reinforcement in graphs. {\it Discrete Appl. Math.},
{\bf 159}(14) (2011), 1443-1446.


\bibitem{hx12}
F.-T. Hu and J.-M. Xu, On the complexity of the bondage and reinforcement problems.
{\it J. Complexity}, {\bf 28}(2) (2012), 192-201.

\bibitem{hx07a}
J. Huang and J.-M. Xu, The total domination and bondage numbers of
extended de bruijn and Kautz digraphs. {\it Comput. Math. Appl.},
{\bf 53}(8) (2007), 1206-1213.

\bibitem{hwx09}
J. Huang, J.-W. Wang and J.-M. Xu, Reinforcement numbers of
digraphs. {\it Discrete Appl. Math.}, {\bf 157}(8) (2009),
1938-1946.

\bibitem{km90}
J. Kok and C.M. Mynhardt, Reinforcement in graphs. {\it Congr.
Numer}, {\bf 79} (1990) 225-231.

\bibitem{kp91}
V.R. Kulli, D.K. Patwari, The total bondage number of a graph, in:
V. R. Kulli (Ed.), {\it Adv. Graph Theory}, Vishwa, Gulbarga,
(1991) 227-235.

\bibitem{mb87}
H. M\"{u}ller, A. Brandst\"{a}dt, The NP-completeness of Steiner tree and dominating set for chordal bipartite graphs.
{\it Theor. Comput. Sci.}, {\bf 53} (1987), 257-265.

\bibitem{plh83}
J. Pfaff, R.C. Laskar, S.T. Hedetniemi, NP-completeness of total and connected domination and irredundance for bipartite graphs.
Technical Report 428, Clemson University, Dept. Math. Sciences, 1983.

\bibitem{ses07a}
N. Sridharan, M.D. Elias, V.S.A. Subramanian, Total bondage
number of a graph. {\it AKCE Int. J. Graphs Combin.}, {\bf 4}(2)
(2007), 203-209.

\bibitem{ses07b}
N. Sridharan, M.D. Elias, V.S.A. Subramanian,
Total reinforcement number of a graph. {\it AKCE Int. J. Graphs
Comb.}, {\bf 4}(2) (2007), 197-202.

\bibitem{xu03}
J.-M. Xu, {\it Theory and Application of Graphs}. Kluwer Academic
Publishers, Dordrecht/Boston/London, 2003.

\bibitem{xu13}
J.-M. Xu, On Bondage Numbers of Graphs: A Survey with Some Comments.
{\it Inter. Journ. Comb.}, {\bf 2013} (2013), Article ID: 595210, 34 pages.

\end{thebibliography}
\end{document}